\documentclass[11pt]{article}
\usepackage{amssymb,latexsym,amsfonts,verbatim,amscd}
\usepackage{amsmath}
\usepackage[all]{xy}

\newtheorem{Def}{Definition}[section]

\newtheorem{Lem}[Def]{Lemma}
\newtheorem{Prop}[Def]{Proposition}
\newtheorem{Cor}[Def]{Corollary}
\newtheorem{Fac}[Def]{Fact}
\newtheorem{Rem}[Def]{Remark}
\newtheorem{Rems}[Def]{Remarks}

\font\nat msbm10 scaled\magstephalf
\def\N{\hbox{\nat\char78}}

\font\mata=msam10 %scaled1095
\def\restr{\mbox{\mata\char22}}
\def\telos{\hfill$\dashv$}

\begin{document}

\title{Sets with dependent elements: A formalization of  Castoriadis' notion of magma}
\author{Athanassios Tzouvaras}

\date{}
\maketitle

\begin{center}
Department  of Mathematics\\  Aristotle University of Thessaloniki \\
541 24 Thessaloniki, Greece \\
e-mail: \verb"tzouvara@math.auth.gr"
\end{center}

\begin{abstract}
We present  a formalization of   collections  that  Cornelius  Castoriadis calls  ``magmas'', especially the property which mainly characterizes them and  distinguishes them from the usual cantorian sets. It is the  property  of their elements to {\em depend}  on other elements, either in a one-way or a two-way manner, so that one cannot occur in a collection without the occurrence of those dependent on it. Such a  dependence  relation on a set $A$ of  atoms (or urelements) can be naturally represented by a pre-order relation $\preccurlyeq$  of $A$ with the extra condition that it contains no minimal elements. Then, working in a mild strengthening of the theory ${\rm ZFA}$, where $A$ is an infinite  set of atoms equipped with a primitive pre-ordering $\preccurlyeq$, the class of magmas  over $A$  is represented  by the class $LO(A,\preccurlyeq)$ of nonempty open subsets of $A$ with respect to the lower topology of $\langle A,\preccurlyeq\rangle$. The non-minimality condition for $\preccurlyeq$ implies that all sets of $LO(A,\preccurlyeq)$ are infinite and none of them is $\subseteq$-minimal. Next the pre-ordering  $\preccurlyeq$ is shifted (by a kind of simulation) to a pre-ordering $\preccurlyeq^+$ on  ${\cal P}(A)$, which turns out to satisfy the same  non-minimality condition as well, and  which, happily,  when restricted to $LO(A,\preccurlyeq)$ coincides with $\subseteq$. This allows us to define a hierarchy $M_\alpha(A)$, along all ordinals  $\alpha\geq 1$, the``magmatic hierarchy'', such that  $M_1(A)=LO(A,\preccurlyeq)$,  $M_{\alpha+1}(A)=LO(M_\alpha(A),\subseteq)$, and  $M_\alpha(A)=\bigcup_{\beta<\alpha}M_\beta(A)$, for a  limit ordinal  $\alpha$. For every $\alpha\geq 1$, $M_\alpha(A)\subseteq V_\alpha(A)$,  where $V_\alpha(A)$ are the levels of the universe $V(A)$ of ${\rm ZFA}$. The class $M(A)=\bigcup_{\alpha\geq 1}M_\alpha(A)$ is the ``magmatic universe above $A$.''  The axioms  of Powerset and  Union (the latter in a restricted version)  turn out to be true in $\langle M(A),\in\rangle$. Besides  it is shown  that three of the five  principles about magmas that  Castoriadis proposed  in his writings (mildly modified and adapted to the needs of formalization), $M2^*$, $M3^*$ and  $M5^*$,  are true of $M(A)$. A selection  of excerpts from these  writings,  in which the concept of magma was first introduced and elaborated, is presented  in the Introduction.
\end{abstract}


\vskip 0.2in


\section{Introduction}
In his seminal book \cite{Ca87} Cornelius Castoriadis\footnote{Cornelius Castoriadis (1922-1997) was a prominent social and political thinker of 20th century. Born in Greece he spent   most of his life in Paris. From 1979 until his death  he was Director of Studies at the \'{E}cole des Hautes \'{E}tudes en Sciences Sociales (EHESS). His monograph \cite{Ca87} is widely considered  his main work. Although  a  humanitarian philosopher by training, he had got  an impressive  solid  background in science, especially in  economics, mathematics and theoretical physics.} is concerned, among many other things, with the way collections of things are presented to (or created by) the human mind. In particular he aims to make us rethink the ``self-evident'' belief of western rationalistic tradition that any  collection of things needs  to be a {\em cantorian} collection, that is a totality of distinct, definite and ontologically independent elements.  His favorite (counter) examples are the  ``totality of meanings'' of a natural language,  the ``totality of one's memories'', and the like. The elements of such  totalities are not entirely  definite, and not fully differentiated  and independent from one another, so one could hardly call them  ``sets'' in the ordinary sense of the word, and thus include them in  the cantorian universe. If, for instance,  $a$ is a particular meaning or memory, one could not fully separate it  from related  meanings or memories, respectively, in the sense that whenever we think of  $a$ as an element of some collection, $a$ inevitably  brings to mind other similar meanings or  memories as members  of the same collection, and as a result  $a$ cannot exist in isolation. A typical  consequence of this is that  the one-element set  $\{a\}$ can hardly make sense, as the cantorian tradition requires. Since such collections abound around us, it would be natural to try to accommodate and comprehend them in the framework of a theory that differs  from that of Cantor.

C. Castoriadis  seems to be  the first thinker who felt the need for some kind of theory that would embrace, even without full rigor, the study of such collections.  He believes that the specific cantorian tradition, which requires collections to comply with the  rules  of standard set theory and rejects  those  that do not as  non-existent,  is rather accidental and due to the early adoption by the western thought of what he calls ``identitary-ensemblistic'' logic (roughly, the  two-valued classical logic  interacting with naive set theory of predicate extensions).

\begin{quotation}
``For the past 25 centuries, Greco-Western thinking has constituted, developed, amplified and refined itself on the basis of this thesis: being is being something determined (einai ti), speaking is saying something determined (ti legein). And, of course, speaking the truth is determining speaking and what is said by the determinations of being or else determining being by the determinations of speaking, and, finally, observing that both are but one and the same. This evolution, instigated by the requirements of one dimension of speaking and amounting to the domination or the autonomization of this dimension, was neither accidental nor inexorable; it corresponded to the institution by the West of thinking as Reason. I call the logic described above identitary logic and also, aware of the anachronism and the stretching of words involved here, set-theoretical logic, for reasons that will soon be apparent. (...)

The logical rudiments of set-theory are important in this respect for, regardless of what may happen in the future from the perspective of mathematics, they condense, clarify, and exemplify in a pure manner what, all the while, was underlying identitary logic, and what, long before this logic was sketched out, constituted an essential and unexpungible dimension of all activity and all social life. These rudiments, indeed, posit and constitute explicitly both the type of logic, in its greatest generality, required by identitary logic and the relations necessary and almost sufficient for this logic to function unhampered and without limit.'' (\cite{Ca87}, pp. 221, 222, 223)
\end{quotation}
Castoriadis believes that  identitary-ensemblistic logic affects greatly our grasping of the reality through the  ``creation'' of  sets out of  something  pre-existent and rather undifferentiated. This  undifferentiated reality out of which the identitary-ensemblistic logic generates sets, classes, objects and properties, is, roughly, what he calls {\em magma}.

\begin{quotation}
``What we seek to understand is the mode of being of what gives itself before identitary or ensemblist logic is imposed; what gives itself in this way in this mode of being, we are calling a {\em magma}. It is obviously not a question of giving a formal definition of it in received language or in any language whatsoever. The following statement, however, may not be unhelpful:

A magma is that from which one can extract (or in which one can construct) an indefinite number of ensemblist organizations but which can never be reconstituted (ideally) by a (finite or infinite) ensemblist composition of these organizations.'' (\cite{Ca87}, p. 343)
\end{quotation}

As mentioned above,  Castoriadis'  favorite  examples of  magmas are  the ``multiplicity  of meanings/significations'' of a natural language and  the ``multiplicity  of one's representations.''
\begin{quotation}
``Let us try then, by means of an accumulation of contradictory metaphors, to give an intuitive description of what we mean by magma (the best intuitive support the reader can present to himself is to think of `all the significations of the English language' or `all the representations of his life'). We have to think of a multiplicity which is not one in the received sense of the term but which we mark out as such, and which is not a multiplicity in the sense that we could actually or virtually enumerate what it `contains' but in which we could mark out in each case terms which are not absolutely jumbled together. (...)

And we have to think of the operations of identitary logic as simultaneous, multiple dissections which transform or actualize these virtual singularities, these components, these terms into distinct and definite elements, solidifying the pre-relation of referral into relation as such, organizing the holding together, the being-in, the being-on, the being-proximate into a system of determined and determining relations (identity, difference, belonging, inclusion), differentiating what they distinguish in this way into `entities' and 'properties' , using this differentiation to constitute 'sets' and `classes.''' (\cite{Ca87}, p. 344)
\end{quotation}

Objects and properties of the  world seem to be the outcome of the ability of our mind for  separation, partitioning  and  individuation. Sets and classes are also products of this very mental mechanism. Castoriadis  mentions time and again  Cantor's well-known ``definition'' of set: ``A set is a collection into  a whole of definite and separate objects of our intuition or our thought. These objects are called `elements' of the set.''  Any specification  of a set is clearly an act of separation and individuation. When we say ``let $X$ be the set of all $x$ such that...,'' we focus on a specific part of the reality, we individualize it and cut it off as a separate object by an act of saying, that is,  a linguistic construct (a formula). On the other hand, the basic  quality of magma, which  differentiates  it from  ordinary sets and classes, is the fact that its ``elements'' are neither fully determined nor fully distinguishable and separable  from one another.

\begin{quotation}
``As a magma, the significations of a language are not the elements of an ensemble subject to determinacy as their mode and their criterion of being. A signification is indefinitely determinable (and the `indefinitely' is obviously essential) without thereby being determined. It can always be marked out, provisionally assigned as an identitary element to an identitary relation with another identitary element (this is the case in designation), and as such be `a something' as the starting point for an open series of successive determinations. These determinations, however, in principle never exhaust it. What is more, they can, and always do, force us to reconsider the initial `something' and lead us to posit it as `something else,' overturning by this very fact, or in order to bring it about, the relations by means of which the initial determination had been made.'' (\cite{Ca87}, p. 346)
\end{quotation}
It is remarkable that a very  similar position is expressed  by John Searle  about mental  states and the content of our consciousness in general:
\begin{quotation}
``One has conscious states such as pains and thoughts only as a part of  living a conscious life, and each state has the identity it has only in relation to other such states. My thought,  for example, about a ski race  I ran long ago, is only that very thought because of its position in a complex network of other thoughts, experiences, and memories. My mental states are internally related to each other in the sense that in order for a mental state to be that state with that character it has to stand in certain relation to the real world.'' (\cite{S98}, p. 42)
\end{quotation}
The above speculations   about magma are clearly  vague. However in \cite{Ca89} Castoriadis  devotes a whole chapter to the subject entitled ``The logic of magmas and the problem of autonomy.'' He starts the chapter with a quotation from a letter of G. Cantor to R. Dedekind: ``Every multiplicity is either an inconsistent multiplicity
or it is a set.'' On this Castoriadis comments:

\begin{quotation}
``To say of a multiplicity that it is inconsistent obviously implies that this multiplicity {\em is,} it is in a certain fashion that remains to be specified and that Cantor does not specify. Clearly, we are not dealing here with an empty set, which is a set in full right, with its place in set theory.

It is toward these inconsistent multiplicities - inconsistent from the standpoint of a logic that claims to be consistent or rigorous - that I turned, starting from the moment, in 1964-1965, when the importance of what I have called the radical imaginary in the human world became apparent to me. Noting that the human psychism cannot be `explained' by biological factors or considered as a logical automaton of no-matter-what richness and complexity. (...)

After various terminological peregrinations - cluster, conglomerate, and others - for this mode of being, as well as the logico-ontological organization it bears, I have ended up with the term {\em magma}. I was later to discover that from 1970 on the editions of Nicolas Bourbaki's {\em Alg\`{e}bre} utilized the term with an acceptation that bears no relation at all to the one I have tried to give it and that is, of course, strictly ensemblistic-identitary in character. As the term, by its connotations, admirably lends itself to what I want to express, and as, dare I say, its utilization by Bourbaki seems to me both rare and superfluous, I have decided to retain it.'' (\cite{Ca89}, pp. 366-368)
\end{quotation}
A little later, after recalling the definition of magma,  Castoriadis says:

\begin{quotation}
``I note in passing that Jean-Pierre Dupuy remarked to me that the `definition'  cited above is unsatisfactory, for it would cover just as well what, to avoid Russell's Paradox, has been called in mathematics a `class.' The objection is formally correct. It does not trouble me much, for I have always thought, and still think, that the `class,' in this acceptation of the word, is a logical artifact constructed {\em ad hoc} to get around Russell's Paradox, and that it succeeds in doing so only by means of an infinite regress. Rather than comment on this `definition,' however, we are going to try here to illuminate other aspects of the idea of magma by exploring the paths (and the impasses) of a more `formal' language. For this, one must introduce a primitive (indefinable and undecomposable) term/relation: the marking (rep\'{e}rer) term/relation, whose valence is at once unary and binary. So, let us suppose that the reader unambiguously understands the expressions: `to mark X;' `X marks Y;' `to mark X in Y' (to mark a dog; the collar marks the dog; to mark or locate the dog in the field). In using this term/relation, I `define' a magma by the following properties:

{\em M1}:  If $M$ is a magma, one can mark, in $M$, an indefinite number
of ensembles.

{\em M2}: If $M$ is a magma, one can mark, in $M$, magmas other than $M$.

{\em M3}: If $M$ is a magma, $M$ cannot be partitioned into magmas.

{\em M4}: If $M$ is a magma, every decomposition of M into ensembles
leaves a magma as residue.

{\em M5}:  What is not a magma is an ensemble or is nothing.'' (\cite{Ca89}, pp. 379-380)
\end{quotation}
Unfortunately, the meaning  of ``marking'' in M1 and  M2  is unclear.
However we guess  that  M1 most likely means: for every magma $M$ there is an indefinite number of  sets $x$ such that $x\subseteq M$. While M2 means: for every magma $M$ there is a magma $N\neq M$ such that $N\subseteq M$. If we interpret M1 this way, then we easily see that M1 and M4 are contradictory. This is pointed out by Castoriadis himself (\cite{Ca89}, p. 383): given a magma $M$, let $X$ be the union of all sets contained in $M$. By M4,  $M\backslash X$ is a magma. But then, by M1, there  is a set $x$ such that $x\subseteq M\backslash X$. This contradicts the fact that $X$ is the union of all sets contained in $M$.

I think that the problem with the principles  M1 and M4 arises from the fact that they both  relate magmas to {\em sets}, which by definition are collections of a different kind, and  in doing so they contradict each other. In contrast M2 and M3 describe  how  magmas relate to other magmas {\em alone}, specifically their submagmas. As for  M5, I construe it as saying: ``What is not a magma is an ensemble and nothing but an ensemble.''  In other words, it suggests that  the classes of magmas and sets exhaust the content of the universe and  are complementary. However in the real world, as well as in that of ${\rm ZFA}$, which  will be used below,  there exist objects/atoms that are non-collections. Therefore  M5 could more realistically be reformulated as follows: ``What is not a magma is a set or an atom.''

Castoriadis thinks  that M3  is the most crucial of the above properties of magmas. He says:
\begin{quotation}
``The third property (M3) is undoubtedly the most decisive. It expresses the impossibility of applying here the schema/operator of separation - and, above all, its irrelevance in this domain. In the magma of my representations, I cannot rigorously separate out those that `refer to my family' from the others. (In other words, in the representations that at first sight `do not refer to my family,' there always originates at least one associative chain that, itself, leads to `my family.' This amounts to saying that a representation is not a `distinct and well-defined being,' but is everything that it brings along with it.) In the significations conveyed by contemporary English, I cannot rigorously separate out those that (not in my representation, but in this tongue [langue] itself) refer in any way at all to mathematics from the others.'' (\cite{Ca89}, p. 381)
\end{quotation}

Let me sum  up. The  description  of magma by Castoriadis through the principles M1-M5 is not sufficiently clear and, most important,  two of these principles, namely M1 and M4, are straightforwardly  contradictory.  Nevertheless  there is an aspect  of the idea that deserves further elaboration. This is the real fact  that we often come across collections, let us call them ``unusual'',  that  differ considerably from  the ``usual'' ones. The latter  include  all   finite collections of things  we deal with in our everyday life, like the collection of students in a class, the collection of books at our bookcase, etc, but also  infinite collections of abstract mathematical entities,  e.g. the collections of integers, rational numbers, real numbers, etc.  The main characteristic of all usual collections is that each one of their elements occurs in absolute (existential) independence from all the rest, so it can be separated and removed from a collection, or added to that, without affecting the occurrence of the others. In contrast,  the  members  of unusual collections come up not in full separation from one another, but rather as  unbreakable  {\em chains} or {\em bunches} of {\em dependent} objects, so that one cannot add or subtract one element without adding or subtracting the  elements depending on it. This dependence is vividly described in the last excerpt above through the example of  things ``that refer to my family'', on the one  hand, and ``those  that do not'', on the other,  and our inability to completely separate one kind from the other.

It is exactly  this type  of {\em collections  with dependent elements}  that we are going  to  consider and formalize  in this paper. As for M1-M5, I propose that, firstly,  M1 and M4 be left out of consideration because of their inconsistency, and, secondly,  the rest  principles  M2, M3, M5 be slightly reformulated as follows:

\vskip 0.2in

{\em M2*}: If $M$ is a magma, there is a magma $N\neq M$ such that $N\subseteq M$.

{\em M3*}: If $M$ is a magma, there is no partition of $M$ into submagmas $M_1$, $M_2$.

{\em M5*}: What is not a magma is a set or an atom.

\vskip 0.1in

\noindent It turns out  that the formalization of magmas that we develop below succeeds in capturing  M2*, M3* and M5*.

\section{Formalizing dependence of  objects in an extension  of ${\rm ZFA}$}
In all Castoriadis's examples of magmas (the collection of meanings of a natural language, the collection of one's mental representations, etc), everyone of their  members shows a clear ``ontological''  dependence on  other members: everyone of them cannot occur in one's mind without the simultaneous occurrence of others. We find this notion of   dependence  interesting and challenging, and it is our purpose in this paper to try to  capture it mathematically. The idea, very roughly,  is to work in the theory ${\rm ZFA}$, which consists of the axioms of ${\rm ZF}$ plus a  set $A$  of non-sets, called ``urelements'' or ``atoms'', which throughout will be referred to for simplicity just as  {\em atoms} (see \cite[p. 250]{Je03} for the formal treatment of this theory). Magmas are going to be represented by  certain subsets of $A$ after the latter will  be endowed with a dependence relation.

Now a first basic  question before going further is whether magmas should be treated as  finite or infinite collections. As we saw in the Introduction, the examples of magmas given by Castoriadis -  the collection of meanings of a natural language, the totality of  memories of a single person and so on -  are all collections that contain elements  of some kind of {\em human resources},  and as such  cannot  be infinite, at least in the strict sense of infinity as we use it in mathematics. That is, we cannot prove e.g. that if $A$ is the collection of memories, then there is an 1-1 function $f:\N\rightarrow A$. But it isn't finite either in the strict sense, because there is no $n\in \N$ for which we can prove that $|A|=n$. The best we can say is that $A$ is {\em potentially infinite}, i.e., a non-infinite collection and yet without specific finite cardinality, unfinished and uncompleted, for which we can potentially discover new elements. Many physical collections of the world are of this kind: the collections of animals, of plants, of species and so on. Now  if for some reason, we want to treat such collections with mathematical means, we cannot treat them as finite, because in such a case we should assign them a specific cardinality $n$. But if $|A|=n$, it would mean that $A$ is a {\em completed} totality, without the possibility to reveal new elements, and this contradicts the previously  described  nature of $A$. Thus necessarily we must treat them as infinite, since in mathematics there is no intermediate state between finite and infinite.\footnote{Notice also that treating physical collections in general as finite sets  may lead to real paradoxes. To give  an example (not of a magma but of a potentially infinite physical collection), let  $X$ be the collection of human beings since their appearance on earth, and let $<$ denote the genealogical order relation on $X$, where  $a<b$ means ``$a$ is the father of the  father of the  father.... of the father of $b$'', where we take finitely many applications of the operation ``father of''. If $X$ were typically finite, then every descending chain $\cdots <a_3<a_2<a_1$ should have a  first element $a_n$. But then $a_n$ either should  be  fatherless, or should have  a father who is not human.  Clearly both conclusions  are absurd. So necessarily $X$ has to be treated as  infinite.}

So we start with an infinite set $A$ of atoms (in ${\rm ZFA}$) and   equip it  with a binary relation which can adequately  capture the most basic properties of dependence, which are just two: reflexivity (every object $a$ depends on $a$) and transitivity (if $a$ depends on $b$ and $b$ depends on $c$, then $a$ depends on $c$). A binary relation with these two properties is a very familiar mathematical object, is called a {\em pre-order relation}, or just a {\em pre-ordering},  and is usually denoted $\preccurlyeq$.\footnote{Of course a more general relation of dependence is sensible, where an object $a$ depends not on  a single element $b$, but rather on a group of elements $\{b_1,\ldots,b_n\}$, but such a relation does not fit to our context.} So we assume that $A$ comes up  with such a relation $\preccurlyeq$. The intended meaning of $a\preccurlyeq b$ is: ``$a$ depends on $b$'',  or  ``$b$ points to $a$,'' or ``$b$ reminds $a$,'' all of which practically mean that every  occurrence of $b$ is followed by the occurrence of $a$.

It follows by the  preceding discussion, that given the relation $\preccurlyeq$ on $A$, ``magmas over $A$'' (with respect to $\preccurlyeq$) are just the collections   $x\subseteq A$ that  are {\em downward closed under $\preccurlyeq$}, i.e., have  the property:
\begin{equation} \label{E:down}
(\forall a,b)(a\preccurlyeq b  \wedge  b\in x \Rightarrow a\in x).
\end{equation}

\begin{Def} \label{D:magma}
{\em Given a pre-ordered set  $\langle A,\preccurlyeq\rangle$, the class $m(A)$ of} magmas over $A$ {\em (with respect to $\preccurlyeq$) consists of  the nonempty subsets of $A$ having property (\ref{E:down}), namely,
$$m(A)=\{x\subseteq A:x\neq \emptyset  \wedge (\forall a,b\in A)(a\in x \wedge b\preccurlyeq a \rightarrow b\in x)\}.$$}
\end{Def}
In my view, the following three  conditions should be met in the treatment of magmas.   Firstly, magmas must coexist together with ordinary sets, as well as with atoms,  in a ``mixed''   universe. (This after all  was explicitly stated as property M5* above.) Secondly, the magmas of the bottom level of the universe must consist exclusively of {\em atoms}, not sets. (Nevertheless,  magmas of higher ranks can be constructed inductively, having as elements magmas of lower ranks.)  Thirdly,  the  dependence relation $\preccurlyeq$ of  atoms should be a {\em primitive} one, that is, not definable from or reducible to other relations of the  ground theory.\footnote{The third condition, concerning  non-definability of $\preccurlyeq$, could be possibly skipped if we assumed that $\preccurlyeq$ can be constructed on $A$ by {\em choice}, but in this case  we should work in ${\rm ZFCA}$ rather than ${\rm ZFA}$. However we think that working in ${\rm ZFA}$, even mildly augmented, is a simpler and more natural option. }

Below we shall treat magmas along the lines of the above  three conditions. So starting with an {\em infinite}  set of atoms $A$ and a primitive pre-ordering $\preccurlyeq$ on it, we shall build the {\em class of  magmas $M(A)$ above $\langle A,\preccurlyeq\rangle$,}  as a subclass of the universe $V(A)$ of the theory ${\rm ZF}$ with atoms,  ${\rm ZFA}$.

Recall that the language of ${\rm ZFA}$ is $L=\{\in, S(\cdot), A(\cdot)\}$, where $S(\cdot)$ and $A(\cdot)$ are the unary predicates (sorts) for sets and atoms, respectively.
\begin{Def} \label{D:universe}
{\em Given the  set $A$ of atoms,} the universe $V(A)$ of ${\rm ZFA}$ {\em is  the class $V(A)=\bigcup_{\alpha\in Ord}V_\alpha(A)$, where:

$V_0(A)=A$,

$V_{\alpha+1}(A)=V_\alpha(A)\cup {\cal P}(V_\alpha(A))$, and

$V_\alpha(A)=\bigcup_{\beta<\alpha}V_\beta(A)$, for limit $\alpha$.}
\end{Def}

\noindent Here, in addition, we need $A$ to carry a  pre-ordering,  so we introduce $\preccurlyeq$  as a new  primitive binary relation symbol, besides $\in$, and extend   $L$ to $L(\preccurlyeq)=L\cup\{\preccurlyeq\}$. As usual, in order to avoid using the sorts $S(\cdot)$ and $A(\cdot)$, we  use   variables  $a,b,c,\ldots$ for atoms, variables $x,y,z,\ldots$ for sets, and also  variables $u$, $v$, $w,\ldots$ that range over both sets and atoms.  The atomic formulas of $L(\preccurlyeq)$  are those of $L$,  plus the formulas  $a\preccurlyeq b$.  For every  $a\in A$,  let $pr(a)$ denote  the set of its  {\em predecessors,}
$$pr(a)=\{b\in A:b\preccurlyeq a\}.$$
In view of the preceding discussion, $pr(a)$ represents the set of elements of $A$ which {\em depend} on $a$. It is well-known that over every pre-ordered set $\langle A,\preccurlyeq\rangle$, the sets $pr(a)$, $a\in A$, form the basis of one of the  natural topologies induced by $\preccurlyeq$, usually called ``lower topology'' for obvious reasons (the corresponding ``upper topology'' has as basis the sets $suc(a)=\{b:a\preccurlyeq b\}$). A set $x\subseteq A$ is said to be {\em open} w.r.t. the lower topology, or {\em lower open},  if for every $a\in x$, $pr(a)\subseteq x$. (It is easy to check that a $x\subseteq A$ is a lower open set if and only if  $A\backslash x$ is an  upper open set.) By the transitivity of $\preccurlyeq$, for every $b\in pr(a)$, $pr(b)\subseteq pr(a)$, so  every $pr(a)$ is open, and we refer to them as ``basic open'' sets, or b.o. sets for short. Let $LO(A,\preccurlyeq)$  denote the set of all {\em nonempty} lower open subsets of $A$. It is obvious that for every family $(x_i)_{i\in I}$ of elements of $LO(A,\preccurlyeq)$,  $\cup_{i\in I}x_i$  belongs to $LO(A,\preccurlyeq)$, and so does  also $\cap_{i\in I}x_i$ whenever it is nonempty.  In particular $A\in LO(A,\preccurlyeq)$. A moment's inspection shows that the sets in $LO(A,\preccurlyeq)$ are exactly those of the class  $m(A)$ defined  in  Definition  \ref{D:magma},  that is,
$$m(A)=LO(A,\preccurlyeq)=\{x\subseteq A: x\neq \emptyset  \wedge (\forall a\in x)(pr(a)\subseteq x)\}.$$
Note that  there can be  $a\neq b$ in $A$ such that $a\preccurlyeq b$ and $b\preccurlyeq a$. This is in accordance with the intuitive meaning of $\preccurlyeq$, and expresses the fact that $a,b$ are {\em mutually dependent}. We  write  then $a\sim b$,  $\sim$ is  an equivalence relation on $A$ and we denote  by $[a]_\sim$, or just $[a]$, the equivalence class of $a$. Obviously, $[a]=[b]$ if and only if  $pr(a)=pr(b)$.

A set $x\in LO(A,\preccurlyeq)$ is said to be {\em $\subseteq$-minimal} or just {\em minimal}, if there is no  $y\in LO(A,\preccurlyeq)$ such that $y\varsubsetneq x$.

\begin{Lem} \label{L:minimal}
Let  $x\in LO(A,\preccurlyeq)$. The following are equivalent.

(i) $x$ is minimal.

(ii) $(\forall a\in x)(x=pr(a))$.

(iii) $(\forall a\in x)(x=[a])$.
\end{Lem}

{\em Proof.} (i) $\Rightarrow$ (ii) Let $x$ be open and minimal and assume that for some $a\in x$, $x\neq pr(a)$. Since, by openness $pr(a)\subseteq x$, it means that $pr(a)\varsubsetneq x$, which contradicts the minimality of $x$.

(ii) $\Rightarrow$ (iii). Assume (ii) is true, and let $a\in x$. Since $[a]\subseteq pr(a)\subseteq x$, it is always true that  $[a]\subseteq x$. For the other inclusion, pick $a\in x$. If $x=\{a\}$, obviously $pr(a)=\{a\}=[a]=x$. If there is $b\in x$ such that $b\neq a$, by (ii) $x=pr(b)=pr(a)$, so $x=[a]=[b]$.

(iii) $\Rightarrow$ (i). We show the contrapositive. Suppose (i) is false, that is  there is an open $y$ such that $y\varsubsetneq x$. Picking $a\in y$, we have $[a]\subseteq pr(a)\subseteq y \varsubsetneq x$, therefore $[a]\neq x$, so (iii) is false.  \telos

\vskip 0.2in

Now for  an arbitrary pre-ordering $\preccurlyeq$ it is clear that  we may have $pr(a)=\{a\}$, for some $a$. In that case $\{a\}$ should be included  in the class of magmas over $A$,  a fact which  is  counter intuitive according to our  previous discussion. But even if $pr(a)\neq \{a\}$ but $pr(a)$ is {\em minimal}, and hence $pr(a)=[a]$, according to Lemma \ref{L:minimal},  we shall have the same problem later,  with respect to the pre-ordering  $\preccurlyeq^+$,  which will be   defined in the next section  on  $LO(A,\preccurlyeq)$. Namely, in that case the singleton  $\{pr(a)\}$ will be open in the lower topology induced by  $\preccurlyeq^+$. So it is necessary to avoid the existence of minimal open sets not only in $LO(A,\preccurlyeq)$, but also in all topologies  induced by the shiftings of $\preccurlyeq$ to  the higher levels of the magmatic hierarchy.

\begin{Lem} \label{L:finite}
If $pr(a)$ is finite, then it contains  minimal open subsets.
\end{Lem}

{\em Proof.}  Let $pr(a)=\{b_1,\ldots,b_n\}$ for some $n\geq 1$. For each $i=1,\ldots,n$, $pr(b_i)\subseteq pr(a)$ is finite, so we can pick one $pr(b_k)$ with the {\em least number} of elements. Then $pr(b_k)$ is minimal. For otherwise $pr(b_k)$ should contain a $b_j$ such that $pr(b_j)\varsubsetneq pr(b_k)$. But then    $|pr(b_j)|<|pr(b_k)|$, a contradiction. \telos

\vskip 0.2in

It follows from the preceding lemma that in order to avoid existence of minimal open sets, it is necessary to impose a condition to $\preccurlyeq$ which implies that every b.o. set $pr(a)$ is infinite.  Given a pre-ordering $\preccurlyeq$ on $A$, let us define  $$a\prec b \ \mbox{iff} \ (a\preccurlyeq b \ \wedge \ b\not\preccurlyeq a) \Leftrightarrow (a\preccurlyeq b \ \wedge \ a\not\sim b).$$
A reasonable  condition which will  guarantee the absence  of minimal open sets in $LO(A,\preccurlyeq)$ is the following:

\vskip 0.1in

(*) \quad \quad $(\forall a\in A)(\exists b\in A)(b\prec a)$.

\begin{Prop} \label{P:suffices}
(i) If $\preccurlyeq$ satisfies (*),  $LO(A,\preccurlyeq)$ does not contain minimal  sets. A fortiori, in view of Lemma \ref{L:finite}, all $pr(a)$ are infinite, and hence all sets in $LO(A,\preccurlyeq)$ are infinite.

(ii) The converse is also true. That is, if (*) fails, then  $LO(A,\preccurlyeq)$ has minimal open sets.
\end{Prop}

{\em Proof.} (i) Suppose $\preccurlyeq$ satisfies (*). It suffices to show that no  b.o. set $pr(a)$ is minimal. Take  a set $pr(a)$. By (*) there is $b$ such that $b\prec a$. Then $b\in pr(a)$, so $pr(b)\subseteq pr(a)$, while  $a\notin pr(b)$, because $a\not\preccurlyeq b$, so $a\in pr(a)\backslash pr(b)$, and  therefore  $pr(b)\varsubsetneq pr(a)$.

(ii) Assume that (*) fails, i.e.  there is $a$ such that $(\forall b)(b\not\prec a)$, or $$(\forall b)(b\preccurlyeq a \rightarrow a\preccurlyeq b).$$
The latter means that $(\forall b)(b\in pr(a)\rightarrow b\in[a])$, or that $pr(a)\subseteq [a]$. Since always $[a]\subseteq pr(a)$, it means that for this  specific $a$, $pr(a)=[a]$. Thus $(\forall c\in pr(a))(pr(a)=[c])$, i.e., condition (iii) of Lemma \ref{L:minimal} holds, and therefore $pr(a)$ is minimal. \telos

\vskip 0.2in

It follows that condition (*) for $\preccurlyeq$ is necessary and sufficient in order for $LO(A,\preccurlyeq)$ not to contain minimal open sets.  It  is exactly  such pre-orderings, satisfying (*),  that   we are  going to deal with below. Moreover this property will be  incorporated in the formal system we shall adopt for the treatment of magmas.  Namely, the formal system in which we work  below differs from  ${\rm ZFA}$  in   the following points:

(a) We  strengthen  the schemes of Separation  and Replacement of  ${\rm ZFA}$ so that they  hold for the  formulas of $L(\preccurlyeq)$ rather than just  $L$. We denote this system  ${\rm ZFA}_\preccurlyeq$. (Notice that without this strengthening, we could not  guarantee, for example,  that the collections  $pr(a)=\{b\in A:b\preccurlyeq a\}$ and $LO(A,\preccurlyeq)$  are sets.)

(b) We add to the axioms of ${\rm ZFA}_\preccurlyeq$ the following statements about $\preccurlyeq$:

\vskip 0.1in

$(D_1)$  $(\forall a)(a\preccurlyeq a)$ (reflexivity).

$(D_2)$  $(\forall a,b,c)(a\preccurlyeq b \wedge b\preccurlyeq c \rightarrow a\preccurlyeq c)$ (transitivity).

$(D_3)$  $(\forall a)(\exists b)(b\prec a)$  (no minimal elements).

\vskip 0.1in

\textbf{Two Remarks about axiom $D_3$.}  1)  $D_3$ suggests  that {\em every} element of $A$ depends on other elements, whereas one may object that, at least in {\em some} magmas,  there may exist {\em independent} elements, i.e., $a\in A$ such that $(\forall b\neq a)(a\not \preccurlyeq b \wedge b\not \preccurlyeq a)$. This is a reasonable objection, provided we agree that such elements are rather exceptions to the rule,  and that the subset of $A$ that consists of  the {\em dependent} elements is still infinite. But then it suffices simply to  take as set of atoms  the set $A'=\{a\in A: (\exists b\neq a)(a\preccurlyeq b \vee b\preccurlyeq a)\}$ instead of $A$, and work with it as before.

2) By condition (*) and Proposition \ref{P:suffices}, axiom $D_3$ is responsible for the  nonexistence of minimal magmas,  and thus, in a sense, the non-trivialization of the magmatic hierarchy. So one might suspect  that $D_3$ is a rather ad hoc axiom. But it is not. Suppose for instance that  $\langle A,\preccurlyeq\rangle$ represents the magmatic collection of  meanings of a natural language and assume that $D_3$ is false for $\preccurlyeq$. Then  $\langle A,\preccurlyeq\rangle$ contains at least one  element $a$ such that   $(\forall b)(b\preccurlyeq\ a \rightarrow a\preccurlyeq b)$. It means  that  $a$ is a linguistic meaning (word) such that  from {\em every} meaning (word) $b$ that descends from $a$, i.e., is generated  from $a$, we recover  $a$. But this is clearly false and unnatural. We can hardly think of a natural language that contains such minimal, that is, {\em terminal} meanings, which  in the long run   do not generate  other meanings that {\em strictly descend} from them. One can easily be convinced about that by remembering that every natural language is a living  organism  constantly changing  and creating new meanings.

\vskip 0.1in

So  the theory we shall be  working in below is
$${\rm ZFA}_D={\rm ZFA}_\preccurlyeq+\{D_1,D_2,D_3\}.$$
In this  theory we shall define the  class of magmas $M(A)$ as a  subclass of the universe $V(A)$. And as the universe of ${\rm ZFA}$ is made of levels $V_\alpha(A)$, for  $\alpha\in Ord$,  the class $M(A)$ will be made also of levels $M_\alpha(A)\subseteq V_\alpha(A)$.  As  first level we  take   the set
\begin{equation} \label{E:firstlevel}
M_1(A)=LO(A,\preccurlyeq),
\end{equation}
so  indeed $M_1(A)\subseteq V_1(A)$.

\begin{Rems} \label{R:useful}
{\rm (i) Since $A\in M_1(A)$, $A$ itself is a magma. If $(x_i)_{i\in I}$ is any family of magmas, then so is $\bigcup_ix_i$,  as well as $\bigcap_ix_i$ if it is nonempty.

(ii) In view of axiom  $D_3$, which is identical to condition (*), and Proposition \ref{P:suffices} (i), all sets of $M_1(A)$ are infinite.

(iii) Every open set $x$ is {\em saturated} with respect to the equivalence relation $\sim$ induced by $\preccurlyeq$. That is, for every $a\in x$ $[a]_\sim\subseteq x$. This follows from  the fact that the b.o. sets $pr(a)$ are  saturated. }
\end{Rems}

Our next step is to define  higher levels $M_\alpha(A)$, for $\alpha\geq 1$, of the magmatic hierarchy. For this task we need to shift the pre-ordering $\preccurlyeq$ of $A$ to the levels ${\cal P}^\alpha(A)$ appropriately.

\section{Shiftings of pre-orderings to powersets}
Given a pre-ordered set $\langle A,\preccurlyeq\rangle$, a  natural way to shift $\preccurlyeq$ to the set ${\cal P}(A)$ is by means of  the relation $\preccurlyeq^+$ of ``simulation'' defined as follows: For  $x,y\subseteq A$, let
\begin{equation} \label{E:shift1}
x \preccurlyeq^+y:\Leftrightarrow (\forall a\in x)(\exists b\in y)(a \preccurlyeq b).
\end{equation}
We call $\preccurlyeq^+$ {\em exponential shifting}, or just {\em shifting},  of $\preccurlyeq$.\footnote{I borrowed the notation $\preccurlyeq^+$  from \cite{Ac88} although, given a relation  $R\subseteq X\times X$, Aczel denotes by $R^+$ the relation on ${\cal P}(X)$ defined by:
$$xR^+y \Leftrightarrow (\forall u\in x)(\exists v\in y)(uRv) \ \& \ (\forall v\in y)(\exists u\in x)(uRv).$$
If in addition  $X$ is a transitive set and  $R\subseteq R^+$, $R^+$   is said to be a ``bisimulation.'' So our definition of $\preccurlyeq^+$ is ``half'' of that of $R^+$. This is because  the above  definition of  $R^+$ is appropriate for symmetric relations $R$, especially equivalences,  while $\preccurlyeq$ is a nonsymmetric relation. }

\begin{Lem} \label{L:easy}
(i) If $\preccurlyeq$ is a  pre-ordering  on $A$,  $\preccurlyeq^+$ is a pre-ordering  on ${\cal P}(A)$. (However, if $\preccurlyeq$ is an order, $\preccurlyeq^+$ need not be so.)

(ii) If $\preccurlyeq$ is a total pre-order, i.e., $a\preccurlyeq b$ or $b\preccurlyeq a$ for all  $a,b\in A$,  so is $\preccurlyeq^+$.
\end{Lem}

{\em Proof.} (i) It is  straightforward that $\preccurlyeq^+$ is reflexive and symmetric. (As a counterexample for the case of orderings, take  $\preccurlyeq$ to be a total order on $A$, and let $x_1$, $x_2$ be distinct cofinal subsets of $(A,\preccurlyeq)$, i.e., $(\forall a\in A)(\exists b\in x_i)(a\preccurlyeq b)$, for $i=1,2$. Then clearly $x_1\preccurlyeq^+ x_2$ and $x_2\preccurlyeq^+ x_1$, while $x_1\neq x_2$.)

(ii) Assume  $\preccurlyeq$ is total and that for $x,y\in {\cal P}(A)$, $x\not\preccurlyeq^+y$. Then $(\exists a\in x)(\forall b\in y)(a\not\preccurlyeq b)$. Since $\preccurlyeq$ is total, it follows  that $(\exists a\in x)(\forall b\in y) (b\preccurlyeq a)$. But then,  by logic alone,   $(\forall b\in y)(\exists a\in x) (b\preccurlyeq a)$, so  $y\preccurlyeq^+x$. \telos

\vskip 0.2in

Since $\preccurlyeq^+$ is a pre-ordering on ${\cal P}(A)$, the sets
$$pr^+(x)=\{y\subseteq A:y\preccurlyeq^+x\}$$
form the b.o. sets of the lower topology $LO({\cal P}(A),\preccurlyeq^+)$. Recall that $LO(A,\preccurlyeq)\subseteq {\cal P}(A)$, so the pre-ordering $\preccurlyeq^+$ applies also to the elements of $LO(A,\preccurlyeq)$. Then the following remarkable relation holds between the sets $pr^+(x)$ and the powersets of $x$.

\begin{Prop} \label{P:connection}
(i) For any $x,y\in {\cal P}(A)$, $y\subseteq x \Rightarrow  y\preccurlyeq^+ x$, so  ${\cal P}(x)\subseteq pr^+(x)$.

(ii) If $x\in LO(A,\preccurlyeq)$, then the converse of (i) holds, i.e., for every $y\in {\cal P}(A)$, $y\preccurlyeq^+ x \Rightarrow  y\subseteq x$, so  $pr^+(x)\subseteq {\cal P}(x)$.

(iii) Therefore for every  $x\in LO(A,\preccurlyeq)$,   $pr^+(x)={\cal P}(x)$.

(iv) In particular, $\preccurlyeq^+\restr LO(A,\preccurlyeq)=\subseteq$, and
$$LO(LO(A,\preccurlyeq), \preccurlyeq^+)=LO(LO(A,\preccurlyeq), \subseteq).$$
\end{Prop}

{\em Proof.} (i). Let $y\subseteq x$ and  pick $a\in y$. Then $a\in x$ and since $a\preccurlyeq a$, it means that there is $b\in x$ such that $a\preccurlyeq b$. Therefore $y\preccurlyeq^+ x$.

(ii) Suppose $x$ is open, let $y\preccurlyeq^+ x$ and pick $a\in y$. Since  $y\preccurlyeq^+x$, there is $b\in x$ such that $a\preccurlyeq b$, i.e., $a\in pr(b)$. But $pr(b)\subseteq x$, since $x$ is open, so $a\in x$ and  therefore $y\subseteq x$.

(iii) and (iv) follow  immediately   from (i) and (ii). \telos

\vskip 0.2in

\begin{Cor} \label{C:minimal}
If $\preccurlyeq$ satisfies $D_3$, i.e., does not contain minimal elements, then so does  $\preccurlyeq^+$ on $LO(A,\preccurlyeq)$.
\end{Cor}

{\em Proof.} Let $x\in LO(A,\preccurlyeq)$. We have to show that there is  $y\in  LO(A,\preccurlyeq)$ such that $y\prec^+ x$, i.e.,  $y\preccurlyeq^+ x$ and $x\not\preccurlyeq^+ y$.  By Proposition  \ref{P:connection} (iv), this amounts to finding  an open $y$ such that    $y\subseteq x$ and $x\not\subseteq y$, or equivalently $y\varsubsetneq x$. Pick an $a\in x$. Then $pr(a)\subseteq x$. By $D_3$, there is a $b$ such that $b\prec a$, hence  $pr(b)\varsubsetneq pr(a)\subseteq x$. Therefore   $pr(b)\varsubsetneq x$ and $pr(b)$ is open.  \telos

\vskip 0.2in

Proposition \ref{P:connection} is crucial for the construction of the levels $M_\alpha(A)$ of the magmatic hierarchy, for every $\alpha\geq 1$. It says that, whatever the starting pre-ordering $\preccurlyeq$ of $A$ is, the restriction of $\preccurlyeq^+$  to  $LO(A,\preccurlyeq)$ is $\subseteq$.  That is, if we set $M_1=LO(A,\preccurlyeq)$,  then $LO(M_1,\preccurlyeq^+)=LO(M_1,\subseteq)$. And for the same reason,  if we set $M_2=LO(M_1,\subseteq)$, the  shifting $\subseteq^+$ of the relation $\subseteq$ of $M_1$ to the sets of $M_2$ is $\subseteq$ again.  And so on with $\subseteq^{++}$, $\subseteq^{+++}$ etc,  for all subsequent levels $M_3$, $M_4,\ldots$, which are defined similarly. As a result,  every such level consists of infinite sets only, ordered by $\subseteq$  with  no minimal element. Moreover, when we reach a limit ordinal $\alpha$, we can take as $M_\alpha$ just the union of all previous  levels, ordered again by $\subseteq$. And then we can proceed further by setting   $M_{\alpha+1}=LO(M_\alpha,\subseteq)$. This way   $M_\alpha$ is defined for every ordinal $\alpha\geq 1$.

\section{The magmatic universe}
Fix a set $A$ of atoms and  a pre-ordering  $\preccurlyeq$ of it which satisfies the axioms of ${\rm ZFA}_D$, in particular $\preccurlyeq$ has no minimal element. In view of Proposition \ref{P:connection} and the remarks at the end of the  last section, we can define the levels $M_\alpha(A)$ of the {\em magmatic hierarchy above $A$}, by induction on $\alpha$ as follows. (For simplicity we write  $M_\alpha$ instead of $M_\alpha(A)$.)

\begin{Def} \label{D:hierarchy}
{\em $M_1=LO(A,\preccurlyeq)$.

$M_{\alpha+1}=LO(M_\alpha,\subseteq)$, for every $\alpha\geq 1$.

$M_\alpha=\bigcup_{1\leq \beta<\alpha}M_\beta$, if $\alpha$  is a limit ordinal.

$M=M(A)=\bigcup_{\alpha\geq 1}M_\alpha$.

\noindent $M$ is said to be the} magmatic universe above $A$ {\em (with respect to the pre-ordering $\preccurlyeq$).}
\end{Def}

Here are some basic  facts about $M_\alpha$'s and $M$.

\begin{Lem} \label{L:basics}
(i) $M_\alpha\subseteq V_\alpha(A)$, $M_1\subseteq {\cal P}(A)\backslash \{\emptyset\}$,  $M_{\alpha+1}\subseteq {\cal P}(M_\alpha)\backslash \{\emptyset\}$ and  $M_\alpha\in M_{\alpha+1}$, for every $\alpha\geq 1$.

(ii) The b.o. sets of the space $M_{\alpha+1}=LO(M_\alpha,\subseteq)$ are the sets
$$pr_\alpha(x)=\{y\in M_\alpha:y\subseteq  x\}={\cal P}(x)\cap M_\alpha.$$
Therefore:
$$x\in M_{\alpha+1} \ \Leftrightarrow \ x\subseteq M_\alpha \wedge x\neq\emptyset \wedge (\forall y\in x)({\cal P}(y)\cap M_\alpha\subseteq x).$$

(iii) The class $M$ is ``almost'' transitive, in the sense that for every $x\in M_\alpha$ for $\alpha\geq 2$, $x\subseteq M$. However for $x\in M_1$,  $x\subseteq A$. So $M\cup A$ is transitive.

(iv) $M$ is a proper subclass of $V(A)$, while $M\cap V=\emptyset$ where $V$ is  the subclass of pure sets of $V(A)$ (that is of $x$'s such that $TC(x)\cap A=\emptyset$).

(v) The inductive definition of $M_\alpha$'s does not reach any fixed point, that is for  every $\alpha\geq 1$, $M_{\alpha+1}\neq M_\alpha$.

(vi) All sets of $M$ are infinite.

(vii) There is no $\subseteq$-minimal set in $M$.
\end{Lem}

{\em Proof.} (i), (ii) and (iii) follow immediately  from the definitions.

(iv) Let $rank(x)$ be the usual rank of sets in $V(A)$, i.e. $rank(x)=\min\{\alpha:x\in V_{\alpha+1}(A)\}$. It is easy to see by induction that for every $\alpha\geq 1$, $rank(M_\alpha)=\alpha$. Concerning the other claim, let  $x\in V\cap M$ be a pure set of least rank $\alpha$. Then  $x\in V_{\alpha+1}\cap M$. If $\alpha\geq 1$,  $x\subseteq V_\alpha\cap M$, and $x\neq \emptyset$. But then  there is   $y\in x$, $y\in V\cap M$ and $rank(y)<rank(x)$, a contradiction. So $rank(x)=0$, which means that $x=\emptyset$ and $x\in V_1\cap M$, a contradiction again since $\emptyset\notin M$.

(v) If  $M_{\alpha+1}=M_\alpha$, then  $M_\alpha=LO(M_\alpha,\subseteq)$, hence  by (i) $M_\alpha\in M_\alpha$, a contradiction. (A bit differently: if $M_{\alpha+1}=M_\alpha$, then $M_\beta=M_\alpha$ for every $\beta>\alpha$, and hence $M=M_\alpha$, so $M$ would be a set, which contradicts (iv).)

(vi) and (vii) follow from the fact that  the initial pre-ordering $\preccurlyeq$ on $A$ does not contain minimal elements (by $D_3$), so by Corollary \ref{C:minimal} so does the relation $\preccurlyeq^+=\subseteq$ on $M_1$, and by induction, so does $\subseteq$ on every $M_\alpha$. Therefore, by Proposition \ref{P:suffices} for every $\alpha\geq 0$,  (a) all sets of $M_{\alpha+1}$ are infinite, and (b) $M_{\alpha+1}$ does not contain $\subseteq$-minimal sets. \telos

\vskip 0.2in

\begin{Rem} \label{notequal}
Notice that as follows from Lemma  \ref{L:basics} (i),  $M_\alpha\subseteq V_\alpha(A)\cap M$. However in general $M_\alpha\neq M\cap V_\alpha(A)$.
\end{Rem}

{\em Proof.} Consider, for example,  the levels $M_1$, $M_2$ of $M$ which are disjoint (see Lemma \ref{L:disjoint} below). Then $M_1\subseteq {\cal P}(A)$ and $M_2\subseteq {\cal P}(M_1)\subseteq {\cal P}^2(A)$. Since ${\cal P}(A)\cup{\cal P}^2(A)\subseteq V_2(A)$, we have  $M_1\cup M_2\subseteq M\cap V_2(A)$. However $M_1\cup M_2\not\subseteq M_2$ because $M_1\cap  M_2=\emptyset$. Therefore $M_2\neq M\cap V_2(A)$. \telos

\vskip 0.1in

By the next three results it is shown that the collection ${\cal P}(x)\cap M$ of submagmas of a magma $x$ is a magma too, that occurs exactly at the next level of that of $x$.

\begin{Lem} \label{L:kalo}
(i) If ${\cal P}(x)\cap M_\alpha\neq \emptyset$, then  ${\cal P}(x)\cap M_\alpha\in M_{\alpha+1}$.

(ii) For every $x\in M$, there is a limit ordinal $\beta$ such that ${\cal P}(x)\cap M={\cal P}(x)\cap M_\beta$.
\end{Lem}

{\em Proof.}  (i) Notice that if $x\in M_\alpha$, then, by Lemma \ref{L:basics} (ii),  ${\cal P}(x)\cap M_\alpha=pr_\alpha(x)$, which is a b.o. set of $LO(M_\alpha,\subseteq)=M_{\alpha+1}$. However we can prove the claim without this assumption. So let $u={\cal P}(x)\cap M_\alpha\neq \emptyset$. Then $u\subseteq M_\alpha$, so  it suffices to show that $u\in LO(M_\alpha,\subseteq)$, i.e., $(\forall z\in u)({\cal P}(z)\cap M_\alpha\subseteq u)$. Pick a $z\in u$. Then $z\subseteq x$, so ${\cal P}(z)\subseteq {\cal P}(x)$, and therefore ${\cal P}(z)\cap M_\alpha\subseteq {\cal P}(x)\cap M_\alpha=u$, as required.

(ii)  $M$ is a (definable) subclass of $V(A)$ and ${\cal P}(x)$ is a set, so by  the Separation scheme of ${\rm ZFA}_D$, ${\cal P}(x)\cap M$ is a set. For each $y\in {\cal P}(x)\cap M$, let $\beta_y=\min\{\gamma\in Ord:y\in M_\gamma\}$. If  $X=\{\beta_y:y\in {\cal P}(x)\cap M\}$,  $X$ is a set of ordinals (by the Replacement axiom), so  $\sup X$ exists, and let $\beta$ be the first limit ordinal such that $\sup X\leq \beta$.  Then $M_\beta=\bigcup_{1\leq \gamma<\beta}M_\gamma$, and  therefore ${\cal P}(x)\cap M\subseteq M_\beta$. But then also ${\cal P}(x)\cap M\subseteq {\cal P}(x)\cap  M_\beta$, so finally ${\cal P}(x)\cap M={\cal P}(x)\cap  M_\beta$ since the reverse inclusion holds  trivially. \telos

\vskip 0.1in

\begin{Prop} \label{P:powrank}
(i) If $x\in M_1$, then ${\cal P}(x)\cap M\subseteq M_1$, and hence ${\cal P}(x)\cap M={\cal P}(x)\cap M_1$. In particular ${\cal P}(A)\cap M=M_1$.

(ii) If $x\in M_{\alpha+1}$, then ${\cal P}(x)\cap M\subseteq M_{\alpha+1}$, for every $\alpha\geq 1$, and hence  ${\cal P}(x)\cap M={\cal P}(x)\cap M_{\alpha+1}$. In particular ${\cal P}(M_\alpha)\cap M= M_{\alpha+1}$.
\end{Prop}

{\em Proof.} (i) Let $x\in M_1$ and pick a $y\in {\cal P}(x)\cap M$. Then $y\subseteq x\subseteq A$.  If $y\notin M_1$, then $y\in M_{\alpha+1}$ for some $\alpha\geq 1$. So $y\subseteq M_\alpha$, and hence $A\cap M_\alpha\neq \emptyset$, which contradicts the fact that $A\cap M=\emptyset$. For the case of $x=A$, we have $A\in M_1$, so ${\cal P}(A)\cap M\subseteq M_1$, but besides  $M_1=LO(A,\preccurlyeq)\subseteq {\cal P}(A)$, therefore ${\cal P}(A)\cap M=M_1$.

(ii) Let $x\in M_{\alpha+1}$, and fix a  $y_0\in {\cal P}(x)\cap M$. We have to show that $y_0\in M_{\alpha+1}$, i.e., $y_0\subseteq M_\alpha$ and $y_0$ is open in $M_\alpha$. Since $y_0\subseteq x\subseteq M_\alpha$, already $y_0\subseteq M_\alpha$. Towards a contradiction, suppose that $y_0$ is not open in $M_\alpha$. It means that the following holds:

(a) $(\exists z\in y_0)({\cal P}(z)\cap M_\alpha\not\subseteq y_0)$.

\noindent Now by Lemma \ref{L:kalo} (ii), there is limit $\beta>\alpha+1$ such that ${\cal P}(x)\cap M={\cal P}(x)\cap M_\beta$, so by Lemma \ref{L:kalo} (i),  ${\cal P}(x)\cap M\in M_{\beta+1}$.  It follows that  $y_0\in M_{\beta+1}$, i.e. $y_0$ is an open subset of $M_\beta$. So

(b) $(\forall  z\in y_0)({\cal P}(z)\cap M_\beta\subseteq y_0)$.

\noindent But since $\beta>\alpha+1$ and $\beta$ is limit, $M_\alpha\subseteq M_\beta$, so

(c) $(\forall z\in y_0)({\cal P}(z)\cap M_\alpha\subseteq {\cal P}(z)\cap M_\beta)$.

\noindent Obviously (b) and (c) contradict (a), and this proves the claim.

In particular for $x=M_\alpha$,  ${\cal P}(M_\alpha)\cap M\subseteq  M_{\alpha+1}$, since $M_\alpha\in M_{\alpha+1}$, but  also $M_{\alpha+1}=LO(M_\alpha,\subseteq)\subseteq {\cal P}(M_\alpha)$, so ${\cal P}(M_\alpha)\cap M= M_{\alpha+1}$. \telos

\begin{Cor} \label{C:PowCon}
(i) For every $\alpha\geq 0$, if $x\in M_{\alpha+1}$ then ${\cal P}(x)\cap M\in M_{\alpha+2}$.

(ii) If $x\in M$ and $x\subseteq M_\alpha$, then $x\in M_{\alpha+1}$  (that is,  every subset of $M_\alpha$ that belongs to $M$, is an open subset of $M_\alpha$).
\end{Cor}

{\em Proof.} (i) By Proposition \ref{P:powrank} (ii), for every $\alpha\geq 0$, if $x\in M_{\alpha+1}$, then ${\cal P}(x)\cap M={\cal P}(x)\cap M_{\alpha+1}$. In addition, by Lemma \ref{L:kalo} (i), ${\cal P}(x)\cap M_{\alpha+1}\in M_{\alpha+2}$.

(ii) It follows again from clause (ii) of Proposition \ref{P:powrank}, in particular from the equality ${\cal P}(M_\alpha)\cap M= M_{\alpha+1}$. \telos

\vskip 0.2in

Proposition \ref{P:powrank} helps us also to compare the  levels of the form  $M_{\alpha+n}$, where  $\alpha$ is a limit ordinal and $n\geq 0$.

\begin{Lem} \label{L:generalL}
For all limit ordinals $\alpha$, $\beta$ such that $\alpha<\beta$, and  all $n\geq 0$, $M_{\alpha+n}\subseteq M_{\beta+n}$.
\end{Lem}

{\em Proof.}  That  $M_\alpha\subseteq M_\beta$ follows from the definition of $M_\alpha$ for limit $\alpha$. We prove the  claim for $n=1$, i.e., $M_{\alpha+1}\subseteq M_{\beta+1}$. Since  $M_{\alpha+1}=LO(M_\alpha,\subseteq)$ and $M_{\beta+1}=LO(M_\beta,\subseteq)$, it suffices to show that every  b.o. subset of $M_\alpha$,  $pr_\alpha(x)$  for some $x\in M_\alpha$, is also a b.o.  subset of $M_\beta$. Pick a $x\in M_\alpha$, so $x\in M_\beta$ too. Since $\alpha$ is a  limit ordinal, $x\in M_{\gamma+1}$ for some  $\gamma<\alpha$, and  by Proposition \ref{P:powrank} (ii), ${\cal P}(x)\cap M={\cal P}(x)\cap M_{\gamma+1}$. Then   $M_{\gamma+1}\subseteq  M_\alpha \subseteq M_\beta$, which implies  that  ${\cal P}(x)\cap M={\cal P}(x)\cap M_\alpha={\cal P}(x)\cap M_\beta$, so
\begin{equation} \label{E:rank}
pr_\alpha(x)={\cal P}(x)\cap M_\alpha={\cal P}(x)\cap M_\beta=pr_\beta(x).
\end{equation}
This proves that every set $pr_\alpha(x)$ is also an open subset of  $M_\beta$, and proves that  $M_{\alpha+1}\subseteq  M_{\beta+1}$. That   $M_{\alpha+n}\subseteq M_{\beta+n}$, for every $n\geq 1$,  is shown by an easy induction and  the help of Proposition \ref{P:powrank} (ii) as before.  \telos

\vskip 0.2in

We give next some results concerning the non-limit levels of $M$. First notice that the levels $M_n$, for finite $n$, have the peculiarity to be pairwise disjoint.

\begin{Lem} \label{L:disjoint}
$M_n\cap M_m=\emptyset$, for  all $n\neq m$.
\end{Lem}

{\em Proof.} Let us set ${\cal P}^*(X)={\cal P}(X)\backslash\{\emptyset\}$, for any set $X$. By Lemma \ref{L:basics} (i), $M_1\subseteq {\cal P}^*(A)$, and  $M_2\subseteq {\cal P}^*(M_1)\subseteq {\cal P}^*({\cal P}^*(A))={\cal P}^{*2}(A)$.  So by induction, for every $n\geq 1$,
\begin{equation} \label{E:one}
M_n\subseteq {\cal P}^{*n}(A).
\end{equation}
On the other hand, $A\cap {\cal P}^{*k}(A)=\emptyset$, for every $k>0$, since $A$ contains non-sets, therefore ${\cal P}^{*n}(A)\cap {\cal P}^{*{n+k}}(A)=\emptyset$, for every  $n$ and every $k>0$, or for all $n\neq m$,
\begin{equation} \label{E:two}
{\cal P}^{*n}(A)\cap {\cal P}^{*m}(A)=\emptyset.
\end{equation}
By (\ref{E:one}) and (\ref{E:two}), $M_n\cap M_m=\emptyset$ for all $n\neq m$.   \telos

\vskip 0.2in

We can give now a uniform  characterization of the sets of  $M_{\alpha+1}=LO(M_\alpha,\subseteq)$, for any limit ordinal $\alpha$.
First observe that if $\alpha$ is limit and $x\subseteq M_\alpha=\bigcup_{1\leq \beta<\alpha}M_\beta$, then $x$ can be written $x=\bigcup_{\beta<\alpha}(x\cap M_{\beta+1})$, because every $y\in x$ belongs to some level of the form $M_{\beta+1}$, for a $\beta<\alpha$.

\begin{Prop} \label{P:partition}
Let $\alpha$ be a limit ordinal. Let $\emptyset\neq x\subseteq M_\alpha$,  $I=\{\beta<\alpha:x\cap M_{\beta+1}\neq \emptyset\}$, and $x_{\beta+1}=x\cap M_{\beta+1}$, for every  $\beta\in I$.  Then $x\in M_{\alpha+1}$ if and only if \begin{equation} \label{E:union}
I\neq \emptyset \ \wedge \  x=\bigcup_{\beta\in I}x_{\beta+1} \ \wedge \ (\forall \beta\in I)(x_{\beta+1}\in M_{\beta+2}).
\end{equation}
\end{Prop}

{\em Proof.} Given  $\emptyset\neq x\subseteq M_\alpha$, clearly $I\neq \emptyset$ and  $x=\bigcup_{\beta\in I}x_{\beta+1}$. So in order to prove (\ref{E:union}) it suffices to prove the equivalence
\begin{equation} \label{E:unionred}
x\in M_{\alpha+1}  \Leftrightarrow  (\forall \beta\in I)(x_{\beta+1}\in M_{\beta+2}).
\end{equation}
$\Rightarrow$ of (\ref{E:unionred}): Assume  $x\in M_{\alpha+1}$, let $\beta\in I$, and pick $y\in x_{\beta+1}$. We have to show that ${\cal P}(y)\cap M_{\beta+1}\subseteq x_{\beta+1}$.  Now $y\in x$ and the assumption $x\in M_{\alpha+1}$ implies  ${\cal P}(y)\cap M_\alpha\subseteq x$. The last inclusion yields ${\cal P}(y)\cap M_\alpha\cap M_{\beta+1}\subseteq x\cap M_{\beta+1}$, which is identical to the required inclusion ${\cal P}(y)\cap M_{\beta+1}\subseteq x_{\beta+1}$.

$\Leftarrow$ of (\ref{E:unionred}): Assume $(\forall \beta\in I)(x_{\beta+1}\in M_{\beta+2})$, where $x_{\beta+1}=x\cap M_{\beta+1}$. It means that for every $\beta\in I$, $(\forall y\in x_{\beta+1})({\cal P}(y)\cap M_{\beta+1}\subseteq x_{\beta+1})$. Pick a $y\in x$. We have to show that ${\cal P}(y)\cap M_\alpha\subseteq x$. Now $y\in x_{\beta+1}$, for some $\beta\in I$, hence $y\in M_{\beta+1}$. By Proposition \ref{P:powrank} (ii), $y\in M_{\beta+1}$ implies that ${\cal P}(y)\cap M\subseteq M_{\beta+1}$, so ${\cal P}(y)\cap M_\alpha={\cal P}(y)\cap M_{\beta+1}$, and since by assumption ${\cal P}(y)\cap M_{\beta+1}\subseteq x_{\beta+1}$, and $x_{\beta+1}\subseteq x$, it follows ${\cal P}(y)\cap M_\alpha\subseteq x$ as required. This completes the proof of (\ref{E:unionred}) and the Proposition. \telos

\begin{Cor} \label{C:follows}
Let $\alpha$ be a  limit ordinal. Then:

(i) For all  $1\leq i\leq n$, $M_i\cap M_{\alpha+n}=\emptyset$. That is, $(\bigcup_{1\leq i\leq n}M_i)\cap M_{\alpha+n}=\emptyset$.

(ii) For  every $n\geq 1$, $M_\alpha\backslash (\bigcup_{1\leq i\leq n}M_i)=\bigcup_{n+1\leq \beta<\alpha}M_\beta\subseteq M_{\alpha+n}$.

(iii) For all $n\geq 0$, $M_{\alpha+n}\not\subseteq M_{\alpha+n+1}$.

\end{Cor}

{\em Proof.} (i)  It suffices to show that given limit $\alpha$, for all $\kappa\geq 1$ and $n\geq 0$, $M_k\cap M_{\alpha+n+k}=\emptyset$. By induction on $k$. Clearly $A\cap M_{\alpha+n}=\emptyset$, so ${\cal P}(A)\cap {\cal P}(M_{\alpha+n})=\emptyset$. Since  $M_1\subseteq {\cal P}(A)$ and $M_{\alpha+n+1}\subseteq {\cal P}(M_{\alpha+n})$ we have $M_1\cap M_{\alpha+n+1}$, so the claim holds for $k=1$. Assume $M_k\cap M_{\alpha+n+k}=\emptyset$. Then $M_{k+1}\subseteq {\cal P}(M_k)$ and $M_{\alpha+n+k+1}\subseteq {\cal P}(M_{\alpha+n+k})$. By the assumption the larger sets in these two relations are disjoint, thus so are the smaller ones.

(ii) The claim  holds for $n=1$, because  every $x\in M_\beta$, for some $2\leq \beta<\alpha$, is a subset of $M_\alpha$ that satisfies condition (\ref{E:union}) of Proposition \ref{P:partition}, so $x\in M_{\alpha+1}$. Then we can continue using induction on $n$. For simplicity, for $n\geq 1$,  assume $M_{n+1}\subseteq M_{\alpha+n}$, and prove that $M_{n+2}\subseteq M_{\alpha+n+1}$.  Let $x\in M_{n+2}$. Then $x\subseteq M_{n+1}$, and $(\forall y\in x)({\cal P}(y)\cap M_{n+1}\subseteq x)$. Now  $y\in x$ implies $y\in M_{n+1}$, so by Proposition \ref{P:powrank} (ii), ${\cal P}(y)\cap M\subseteq M_{n+1}$. By this fact and  the induction assumption $x\subseteq M_{n+1}\subseteq M_{\alpha+n}$, it follows that $x\in M_{\alpha+n+1}$.

(iii) By (i) and (ii) above, $M_\alpha\not\subseteq M_{\alpha+1}$ because $M_1\subseteq M_\alpha$ while $M_1\cap M_{\alpha+1}=\emptyset$, and $M_{\alpha+n}\not\subseteq M_{\alpha+n+1}$ because $M_{n+1}\subseteq M_{\alpha+n}$, while $M_{n+1}\cap  M_{\alpha+n+1}=\emptyset$. \telos

\vskip 0.2in

It follows from  clause (iii) of the preceding corollary that,  in contrast to the limit levels $M_\alpha$,  which, by definition, contain  the elements of all  lower levels,  the successor levels do not behave  cumulatively. Each level $M_{\alpha+n}$, for a limit $\alpha$ and $n\geq 1$, always  omits elements of lower levels.
One may infer from the comparison given above that $M_{\alpha+n}$ and $M_{\alpha+n+1}$ differ only with respect to elements of the initial levels $M_n$ of the hierarchy, i.e. elements of $M_\omega$.  But this is not true. For example pick a  $x\in M_2$ and  a $y\in  M_3$.  Then clearly $x\cup y\in M_{\omega+1}$, according to the characterization given in Proposition \ref{P:partition}, while $x\cup y\notin M_\omega$, so $x\cup y\in M_{\omega+1}\backslash M_\omega$.   But also $x\cup y\notin M_{\omega+2}$ either, because otherwise $x\cup y\subseteq  M_{\omega+1}$. Since $x\subseteq M_1$,  that would mean  that   $M_1\cap M_{\omega+1}\neq\emptyset$, which  is false according to clause (i) of the preceding Corollary.  This shows that  $M_{\omega+1}\backslash M_\omega\not\subseteq M_{\omega+2}$.

\vskip 0.1in

We turn now to another more standard point of view from which the class $M$ could be looked at: the point of view from which $M$ is seen as  an $\in$-structure and so questions are raised as to what set-theoretic properties  of the language $L_0=\{\in\}$ could be satisfied in $\langle M,\in\rangle$. Such a question about the truth of some of the axioms of ${\rm ZF}$ in $\langle M,\in\rangle$  is reasonable. There is however a technical difficulty with sentences of $L_0$ in $M$, because of  the  lack of transitivity due to the  bottom level $M_1$, since for every $x\in M_1$, $x\cap M=\emptyset$. In view of this, given $x,y\in M_1$, the truth of simple properties like $x=y$ and $x\subseteq y$  cannot be  expressed inside $M$ by the usual formulas.  The problem is fixed if we add to $M$ the atoms and work  in  $M^*=M\cup A$ rather than $M$,  with   language $L=\{\in, S(\cdot), A(\cdot)\}$, or with sorted variables.

Still, of course,  we do not   expect  $M^*$ to satisfy many of the closure properties expressed through  the axioms of ${\rm ZFA}$. Of these axioms Extensionality and Foundation do hold in  $\langle M^*,\in\rangle$  because the latter is a transitive substructure  of $\langle V(A),\in\rangle$. But since by construction  $\emptyset\notin M$,  the Emptyset axiom fails. So  does the Pairing axiom because every $x\in M$ is an infinite set. The Infinity axiom fails too, since  the ``measure'' by which we decide   infinity of a set is $\omega=\{0,1,\ldots\}$, and this  is not a resident of $M^*$. Finally the truth of Separation and  Replacement in $M^*$ is obviously out of the question. Nevertheless, the rest two axioms,  Powerset  ($Pow$) and  Union ($Un$), are indeed  true (the second one not quite). We begin  with $Pow$ which is a direct consequence of a result established previously.

\begin{Prop} \label{P:pow}
$M^*\models Pow$. Intuitively, the collection of submagmas of a magma is again a magma (of the next higher rank.)
\end{Prop}

{\em Proof.} We have to show that $M^*\models (\forall x)(\exists y)(y={\cal P}(x))$, or $(\forall x\in M)(\exists y\in M)(y={\cal P}(x)\cap M)$. However  this follows immediately from Corollary \ref{C:PowCon} (i): if  $x\in M_{\alpha+1}$, then   ${\cal P}(x)\cap M\in M_{\alpha+2}$.  \telos

\vskip 0.2in

In contrast to the powerset operation, which ``goes upward'', the union operation ``goes downward'' and may lead out of $M$ if $\cup x$ hits the bottom level that consists of atoms.  For example if $x\in M_1$ then $x\subseteq A$, so  according to the formal definition of $\cup x$, $\cup x=\emptyset \notin M$. Perhaps one may guess that this concerns the elements of $M_1$ only and that  $(\forall x\in M\backslash M_1)(\cup x\in M)$. But still this is not true. For example, as follows  from Proposition \ref{P:partition}, $M_1\cup M_2\in M_{\omega+1}$  while $\cup (M_1\cup M_2)=(\cup M_1)\cup (\cup M_2)=A\cup M_1$. This   does not belong to $M$ either, because on the one hand clearly  $A\cup M_1\notin M_1$, and on the other if we assume $A\cup M_1\in M_{\alpha+1}$, for some $\alpha\geq 1$,  then $A\cup M_1\subseteq M_\alpha\subseteq M$, which is false because  $A\cap M=\emptyset$.

In fact, since the elements of every $x\notin M_1$ are open sets, and every union of open sets of the {\em same space} $LO(M_\alpha,\subseteq)$ is open again, it follows that if $x\subseteq LO(A,\preccurlyeq)$, or $x\subseteq LO(M_\alpha,\subseteq)$, for some $\alpha\geq 1$, then $\cup x\in LO(A,\preccurlyeq)$, or $\cup x\in LO(M_\alpha,\subseteq)$, respectively. This situation occurs exactly  when  $x\in M_{\alpha+2}$,  for some $\alpha\geq 0$, and  gives a sufficient condition in order for $\cup x$ to belong to $M$. It turns out that this condition is also necessary.  In fact the next proposition describes precisely, with two equivalent conditions,  the elements of $M$ whose unions belong to $M$.

\begin{Prop} \label{P:guaran}
Let $x\in M$. The following conditions are equivalent.

(i) $x\in M_{\alpha+2}$, for some $\alpha\geq 0$.

(ii) $\cup x\in M$.

(iii) $x\subseteq M_1$ or $x\cap M_1=\emptyset$.
\end{Prop}

{\em Proof.} (i)$\Rightarrow$(ii)  Assume first that $\alpha\neq 0$ and $x\in M_{\alpha+2}$. Then $x\subseteq M_{\alpha+1}=LO(M_\alpha,\subseteq)$, so $x$ is a family of open subsets of $M_\alpha$. The union $\cup x$ of this family is open too, so  $\cup x\in M_{\alpha+1}$.   If $\alpha=0$, then  $x\in M_2$, or  $x\subseteq M_1=LO(A,\preccurlyeq)$, so  $\cup x\in M_1$. In both cases $\cup x\in M$.

(ii)$\Rightarrow$ (iii): We prove the contrapositive. Suppose $x\not\subseteq M_1$ and $x\cap M_1\neq \emptyset$. If $x_1=x\cap M_1$, then $x=x_1\cup x_2$, where $x_1,x_2\neq \emptyset$, $x_1\subseteq M_1$ and $x_2\subseteq M\backslash M_1$. Therefore  $\cup x=(\cup x_1)\cup (\cup x_2)$, where  $\cup x_1\subseteq \cup M_1=A$ and  $\cup x_2 \cap A=\emptyset$. So if $\cup x\in M$, necessarily  $\cup x\subseteq M_\alpha$ for some $\alpha \geq 1$, and hence $\cup x_1\subseteq M_\alpha$. But then  $A\cap M_\alpha\neq\emptyset$, a contradiction.

(iii)$\Rightarrow$ (i): Let $x\subseteq M_1$. By Corollary \ref{C:PowCon} (ii), $x\in M$ and $x\subseteq M_1$ imply $x\in M_2$, so we are done. Let now $x\cap M_1=\emptyset$. Then $x\subseteq M\backslash M_1$, and if $x\in M_{\alpha+1}$, then $x\subseteq M_\alpha\backslash M_1$. Without loss of generality we may take $\alpha$ to be a limit ordinal. Then by Corollary \ref{C:follows} (ii), $M_\alpha\backslash M_1\subseteq M_{\alpha+1}$, so $x\subseteq M_{\alpha+1}$. By Corollary \ref{C:PowCon} (ii) again, $x\in M_{\alpha+2}$. \telos

\vskip 0.2in

On the other hand every level of the form $M_{\alpha+1}$, for limit $\alpha$, contains sets $x$ such that $\cup x\notin M$.

\begin{Fac}
If   $\alpha=0$, or $\alpha$ is a limit ordinal,  there exists $x\in M_{\alpha+1}$ such that $\cup x\notin M$.
\end{Fac}

{\em Proof.} We saw above  that for every $x\in M_1$, $\cup x\notin M$. Also if, for example, $x=M_1\cup M_2$, then $x\in M_{\omega+1}$ and $\cup x\notin M$. By Lemma \ref{L:generalL} $M_{\omega+1}\subseteq M_{\alpha+1}$, for every limit $\alpha$, so $M_1\cup M_2\in M_{\alpha+1}$ too. \telos

\vskip 0.2in

We close  here the description of the technical features of the magmatic universe and  come to  the question about the extent to which   the class $M$  actually satisfies  some, or all, of   Castoriadis' intuitive principles M1-M5,  which he proposed  as main  properties of magmas. Recall however (see Introduction) that we decided to leave out M1 and M4 as inconsistent, and reformulate  slightly the rest of them into  M2*, M3*, M5*. Given this adjustment,  the answer to the question is yes, M2*, (a weak form of) M3* and M5*  are true in $M$. In our formalization some magmas can be called {\em basic}, if they generate all the rest of the same level. They are just the b.o. sets $pr(a)\in M_1$, for $a\in A$, and  $pr_\alpha(x)\in M_{\alpha+1}$, for $x\in M_\alpha$ and  $\alpha\geq 1$.

\begin{Prop} \label{P:hold}
The following statements  are true about $M$:

(i) (M2*): For every magma $x$ there is a magma $y\neq x$ such that $y\subseteq x$.

(ii) (weak M3*): If $x$ is a basic magma,  there is no finite partitioning of $x$ into submagmas.

(iii) (M5*): What is not a magma is a set or an atom.
\end{Prop}

{\em Proof.} (i) Let $x\in M$ be a magma, and let $x\in M_{\alpha+1}=LO(M_\alpha, \subseteq)$, for some $\alpha\geq 1$, or  $x\in M_1=LO(A,\preccurlyeq)$.  We know (see Lemma \ref{L:basics})  that none  of  these topologies  contains minimal open sets, so there is always a $y\in M_{\alpha+1}$ such that $y\varsubsetneq x$. Such a $y$ is a proper submagma of $x$.

(ii) Given a basic  magma $pr_\alpha(x)$, assume for simplicity that $pr_\alpha(x)=y_1\cup y_2$, where $y_1$, $y_2$ are disjoint submagmas, that is  disjoint  open  sets  $y_1$, $y_2$.  But then  either $x\in y_1$ or $x\in y_2$, hence either $pr_\alpha(x)\subseteq y_1$, or $pr_\alpha(x)\subseteq y_2$, both  of which contradict the fact that $pr_\alpha(x)=y_1\cup y_2$ and $y_1\cap y_2=\emptyset$.

(iii) This claim follows simply from the fact that the class $M=M(A)$ of magmas above $A$ is a subclass of the set-theoretic universe  with atoms $V(A)$. So if $u\notin M(A)$, then necessarily  either $u$ is a classical set of $V(A)$, or $u\in A$. \telos

\vskip 0.2in

{\textbf{Conclusions and some possibilities for future work.} A reviewer wrote that ``Castoriadis himself would probably reject this attempt of formalizing magmas. In particular, he seems to think that set theory might not be the right framework for such a formalization''.  I largely agree. The above presented  formalization in no way should  be  construed as an attempt to capture or interpret mathematically  Castoriadis' thoughts and intuitions in some ``authentic'' fashion; not even to capture his theory of magmas {\em in all of its aspects}. After all I think  Castoriadis did not keep secret his low appreciation for the pair ``classical logic/classical set theory'', to which he, disdainfully enough, referred throughout his writings (among them in some  excerpts cited in the Introduction) as  ``identitary-ensemblistic logic'' -- as if there were plenty  of {\em superior} logics and theories of collections around.

The concept of magma was of interest to me simply because it  challenged the fundamental  assumption of classical set theory that every object  behaves  {\em independently }  of all the rest with respect to $\in$. And while this is unquestionable  for {\em pure sets}, it is not that  obvious for  atoms (urelements) and perhaps for sets containing  atoms. So my interest in magmas stemmed primarily  from the relation of {\em dependence} that underlies their elements, and the challenge to build a structure for an inclusive class of sets with dependent elements. Towards this aim the only {\em secure} tools available to me were just the old good ``identitary-ensemblistic  logic''. Nevertheless  I would be very much interested to see other, alternative approaches to the concept, through different  logical and/or set-theoretical environments.

Concerning the possibilities for future work I can see two directions:

(a)  To try to construct in $M$ ``magmatic analogs'', or ``magmatic codes''  of some standard set-theoretical objects, like ordered pair, function, natural number, $\N$, etc. The starting point could be the observation that the magmatic analogs of the simplest sets  $\{a\}$ of ${\rm ZFA}$, for $a\in A$, are the basic open sets $pr(a)$ of $M$ ($pr(a)$ is the smallest magma that contains $a$). More generally the analog of the $n$-element set $\{a_1,\ldots,a_n\}$ is the magma $pr(a_1)\cup\cdots\cup pr(a_n)$. I cannot say how far this magmatic coding can go and how complicated and intuitively attractive would be. I just wonder   how our mathematical intuition would be if we were living not in the present ``separative'' world  of independent entities, but  in a ``non-separative'' world of dependent entities and  ``thick'' magmatic objects.

(b) The other possible direction  is the {\em syntactic} approach to magmas, i.e., an attempt to axiomatize their basic properties by  a  set of axioms $T$ in a suitable language and then look for models of $T$. Normally the language must include, except  $\in$, a primitive binary relation $D(x,y)$ for dependence. However to  this direction, of  help might be an old paper  of H. Skala, ``An alternative way of avoiding set-theoretical paradoxes'',  {\em Zeitsch. f. math. Logik und Grundlagen d. Math.}  20 (1974), pp. 233-237 (the former title of {\em Mathematical Logic Quarterly}), which, strangely enough,  is related to  sets with dependent elements. This is  because Skala's  axioms, which overlap with those of ${\rm ZF}$, without forming  either a subsystem or an extension of ${\rm ZF}$,  do not include Pairing,
and leave  room for  objects $a\neq  b$ such that $(\forall x)(a\in x\rightarrow b\in x)$. If we denote the last relation between $a$ and $b$ by $b\preccurlyeq a$, then $\preccurlyeq$ is a pre-ordering, and hence a dependence relation.

\vskip 0.2in

{\textbf{Acknowledgements.} I would like to thank Georgios Evangelopoulos for encouraging me to reconsider and rework  an old, unfinished draft  on the subject that I had left aside for years. The outcome of that reworking is the present article. I would like also to thank  the two anonymous reviewers  for constructive criticism and valuable remarks and suggestions.

\end{document}